\documentclass[12pt]{amsart}
\usepackage{amssymb,amsmath}
\numberwithin{equation}{section}
    \def\u+{u^{\tau\,+}}
\def\di{\partial}
\def\dib{\bar\partial}
\numberwithin{equation}{section}
\def\upiuh{u^{+\,(h)}}
\def\upiu0{u^{+\,(0)}}
\def\simleq{\underset\sim<}
\def\simgeq{\underset\sim>}
\def\simle{\underset\sim<}

\def\T{\text}
\def\1#1{\overline{#1}}
\def\2#1{\widetilde{#1}}
\def\3#1{\widehat{#1}}
\def\4#1{\mathbb{#1}}
\def\5#1{\frak{#1}}
\def\6#1{{\mathcal{#1}}}
\def\C{{\4C}}
\def\R{{\4R}}

\def\B{\Bbb B}

\def\Re{{\sf Re}\,}

\def\phi{\varphi}

\def\Om{\Omega}

\newtheorem{Thm}{Theorem}[section]
\newtheorem{Cor}[Thm]{Corollary}
\newtheorem{Pro}[Thm]{Proposition}
\newtheorem{Lem}[Thm]{Lemma}
\theoremstyle{definition}\newtheorem{Def}[Thm]{Definition}
\theoremstyle{remark}
\newtheorem{Rem}[Thm]{Remark}
\newtheorem{Exa}[Thm]{Example}

\def\Label#1{\label{#1}}
\def\bl{\begin{Lem}}
\def\el{\end{Lem}}
\def\bp{\begin{Pro}}
\def\ep{\end{Pro}}
\def\bt{\begin{Thm}}
\def\et{\end{Thm}}
\def\bc{\begin{Cor}}
\def\ec{\end{Cor}}
\def\bd{\begin{Def}}
\def\ed{\end{Def}}
\def\br{\begin{Rem}}
\def\er{\end{Rem}}
\def\be{\begin{Exa}}
\def\ee{\end{Exa}}
\def\bpf{\begin{proof}}
\def\epf{\end{proof}}
\def\ben{\begin{enumerate}}
\def\een{\end{enumerate}}

\def\1alpha{[\frac1\alpha]}
\def\T{\text}
\def\R{{\Bbb R}}

\def\C{{\Bbb C}}

\numberwithin{equation}{section}
\def\T{\text}

\newcommand{\no}[1]{\|{#1}\|}
\newcommand{\NO}[1]{\|{#1}\|^2}

\newtheorem{theorem}{Theorem  }[section]
\newtheorem{definition}[theorem]{Definition }
\newtheorem{lemma}[theorem]{Lemma  }
\newtheorem{proposition}[theorem]{Proposition  }
\newtheorem{corollary}[theorem]{Corollary }
\newtheorem{example}[theorem]{\it Example }

\begin{document}
\title[The Hardy-Littlewood lemma...]{The Hardy-Littlewood Lemma and the estimate of the $\dib$-Neumann problem in a general
norm } 
\author[ S.~Pinton  ]
{Stefano Pinton }
\address{Dipartimento di Matematica, Universit\`a di Padova, via 
Trieste 63, 35121 Padova, Italy}
\email{pinton@math.unipd.it}
\maketitle
\begin{abstract}
We prove a generalized Hardy-Littlewood lemma on a non-smooth domain in ``$f$-norm" and give an application to a corresponding estimate for the $\dib$-Neumann problem by means of suitable weights.

\noindent
MSC: 32F10, 32F20, 32N15, 32T25 
\end{abstract}

\section{Introduction}
\Label{s0}
For a general, smooth function $f=f(t),\,\,t\in\R^+$, with $f$ increasing and a  domain $\Om\subset \R^{2n}$, defined by $r<0$, we 
take ``adapted" coordinates $(a,r)\in\R^{2n-1}\times\R$ and consider the tangential pseudodifferential operator $f(\Lambda)$ with symbol $f(\Lambda_\xi)$ where $\Lambda_\xi$ is the Bessel potential in the variables $\xi$'s dual to $a$'s. 
When $b\Om\in H^{2+n +\epsilon}$, we prove
\begin{equation}
\begin{cases}
\Label{1.1}
\no{f(\Lambda)u}^\Om\le c \no{\delta f(\delta^{-1}) \nabla u}^\Om+\underset{j=1}{\overset n\sum}\no{\bar L_j u}^\Om\quad \T{if $\frac{f}t$ is decreasing},
\\
\no{f(\Lambda)u}^\Om\le c \no{ f(\delta^{-1})  u}^\Om+\underset{j=1}{\overset n\sum}\no{\bar L_j u}^\Om\quad\T{if $\frac{f}{t^\eta}$ is decreasing and $\eta<\frac12$},
\end{cases}
\end{equation}
where $\bar L_j,\,\,j=1,...,n$ is a basis of $ (0,1)$   vector fields and $\delta=\delta(z)$ is the distance function to $b\Om$.
This is similar to the classical  Hardy-Littlewood Theorem in which  $f$ is a fractional $s$-power and hence $\no{ u}_f$ is the Sobolev norm $||u||_s$; note that in the classical statement one has to replace $\sum_j\no{\bar L_j u}^\Om$ by $\no{u}^\Om$ on the right hand side. For the general theory of fractional Sobolev spaces $H^s$, we refer to Lions-Magenes \cite{LM72} for $b\Om$ smooth and Jerison-Kenig \cite{JK95} for $b\Om$ Lipschitz.
The comparison of the general $f$ norm with the usual Sobolev norm is in order. Now, the estimate for the Sobolev norm $\no{u}^\Om_s$ is obtained by proving its equivalence with a sort of trace norm (\cite{LM72} formula (10,5) and \cite{CS00} Appendix C formula (3.4)). This is a consequence of the weighted Hardy inequality \cite{LM72} (10,9). We point out that when the weight is no longer a fractional power but a more general function such as for instance a logarithm, then the inequality breaks down (\cite{T69}). Therefore, for a general $f$, the  appropriate method seems to be that of straightening the boundary and taking trace and harmonic extension.

 We then pass to consider in Section~\ref{s2} the $\dib$-Neumann problem on a bounded pseudoconvex domain $\Om$ with $H^{2+n +\epsilon}$ boundary. We consider  the problem whether the $f$-property, that is, the existence of a uniformly bounded family of weights $\{\phi_\epsilon\}$ whose Levi form satisfies $\di\dib \phi_\epsilon\simgeq f^2(\epsilon^{-1})$ on the $\epsilon$-strip of $\Om$ about $b\Om$, implies  the ${ f_1}$-estimate for $f_1:=\frac f{\log f}$, that is, $\no{ u}_{f_1}\simleq \no{\dib u}+\no{\dib^* u}$ for any form $u$ in the domain of $\dib^*$. When $b\Om\in C^\infty$ this is proved by Catlin \cite{C87} for $f(\epsilon^{-1})=\epsilon^{-\eta}$ and by Khanh \cite{Kh09} for a general $f$.
We get a new proof of this result which has also the advantage of taking a natural generalization to a  boundary which is not necessarily $C^\infty$ but just $H^{2+n +\epsilon}$. First, out of the family $\{\phi_\epsilon\}$, we construct a single weight $\phi$ such that $\di\dib \phi\simgeq f_1(\delta^{-1})$ (it is here that a loss from $f$ to $f_1$ occurs). We then use the basic estimate in the domain $\Om$ weighted by $\phi$ satisfying $\di\dib\phi\simgeq f_1(\delta^{-1})$ and get
\begin{equation}
\Label{1.2}
\no{f_1(\delta^{-1})u}^{\Om}+\sum_{j=1}^n\no{\bar L_j u}^{\Om}\simleq \no{\dib u}^{\Om}+\no{\dib^*u}^{\Om}.
\end{equation}
When $\frac{f}{t^\eta},\,\,\eta<\frac12$, if we  plug \eqref{1.2} with the second of \eqref{1.1}, in which the constant $c$ is uniform with respect to $m$,  we get the $f_1$-estimate for $u$ on $\Om$. 

\vskip0.3cm
\noindent
{\it Aknowledgements}
The paper was accomplished in May 2012 while the  author was supported by a INDAM-GNAMPA grant.

\section{The Hardy-Littlewood Lemma in general norm.}
\Label{s2}
We identify $\C^n$ to $\R^{2n}$ by $z=x+iy\to (x,y)$, take a domain $\Om\subset\C^n$ defined by $r<0$, and complete 
$r$ to a full system of local coordinates $(a,r)$; all through the paper the regularity of the coordinates is at least $C^2$. Let $\xi$ be dual coordinates to the $a$'s, and denoe by $\Lambda_\xi=(1+|\xi|^2)^{\frac12}$ for $|\xi|^2=\sum_j\xi_j^2$, be the standard elliptic symbol of order 1, $f(\Lambda_\xi)$ a  symbol obtained by composition with a general increasing function $f$ and $f(\Lambda)$ the associated tangential pseudodifferential operator defined by
$$
f(\Lambda)u(a,r):=(2\pi)^{-(2n-1)}\int_{\R^{2n-1}}e^{i\langle a,\xi\rangle}f(\Lambda_\xi)\tilde u(\xi,r))d\xi
$$
 where $\tilde u$ is the partial Fourier transform in $\R^{2n-1}$. 

  We also consider  some more general operator $\Psi$ 
with symbol $\psi=\psi(x,r;\xi)$. This is for instance the operator of microlocal decomposition of $u$ into $u^++u^-+u^0=\zeta \Psi^+u+\zeta\Psi^-u+\zeta
\Psi^0u$ for $\zeta\in C^\infty_c(U)$ with $\zeta\equiv1$ in a neighborhood of $\T{supp}\,u$, associated to a conical partition 
of the unity $\psi^+(\xi)+\psi^-(\xi)+\psi^0(\xi)\equiv1$ with $\T{supp}\,\Psi^\pm\subset\{\xi: \pm\xi_{2^n-1}>\frac12|\xi|\}
$. (Here $\di_{x_{2n-1}}=J\di_r|_{z_o}$ for $J$ representing the complex structure in $\C^n$.)
The operation of taking symbols from operators is denoted by $\sigma$ and its reverse by $Op$.
We recall that for the composition we have the rule $\sigma(\Psi_1\circ\Psi_2)=\psi_1\cdot\psi_2-\sum_{|\alpha|>1} \di_x^\alpha\psi_1\di_\xi^\alpha\psi_2$. In particular, $\sigma(f(\Lambda)\Lambda^{-1})=f(\Lambda_\xi)\Lambda_\xi^{-1}$,
but, because of the low regularity of our boundary, we will be not able to apply this rule to the tangential derivatives to the boundary.
 Related to these, is the harmonic extension operator on positive microlocalization that we denote by $u^{+\,(h)}$. To define it, let $L_j,\,\,j=1,...,n$ be a basis of $T^{1,0}\C^n$ with $L_jr=\kappa_{jn}$ (the Kronecker symbol); thus $L_j|_{b\Om},\,\,j\le n-1$ is a basis for $T^{1,0}b\Om$. Let $T:=\frac1{2i}(L_n-\bar L_n)$ and $\di_r:=\frac12(L_n+\bar L_n)$ be the totally real tangential and normal vector fields respectively: they are related by $\di_r=\bar L_n+iT$. We define
\begin{equation}
\Label{2.1}
u^{+\,(h)}:=(2\pi)^{-2n+1}\int_{\R^{2n-1}}e^{ix\xi}e^{r\sigma(\dot T)}\tilde u(\xi,0)d\xi,\quad r<0,
\end{equation}
where $\sigma(\dot T)=\sigma(T)(x,0;\xi)$. 
This can be viewed as a tangential pseudodifferential operator with parameter $r$.
Strictly speaking, the harmonic extension is an operation from $b\Om$ to $\Om$; thus, in the literature it is more usual the notation $(\dot u)^{+\,(h)}$ where ``dot" is the operation of restriction to $b\Om$.

We now wish to compare the action of $f(\Lambda)$ with that of $\delta f(\delta^{-1})\nabla$  and $f(\delta^{-1})$. 
In this discussion, $f$ is a smooth increasing function, $u$ is a $C^\infty(\bar\Om)$ function with support in a neighbourhood $U$ of a boundary point $z_o\in b\Om$ and $b\Om$ is $H^{2+n+\epsilon}$.
Before stating our results, we need some preparation about the action of $f(\Lambda)$ over $\upiuh$.
\bl
\Label{l10.1}
We have
\begin{equation}
\begin{cases}
\Label{10.1}
\int_{-\infty}^0r^2f^2(\frac t{-r})e^{2r}dr\sim f^2(t)\quad\T{if $\frac{f}t$ is decreasing},
\\
\int_{-\infty}^0f^2(\frac t{-r})e^{2r}dr\sim f^2(t)\quad\T{if $\frac{f}{t^\eta}$ is decreasing and $\eta<\frac12$}.
\end{cases}
\end{equation}
\el
\bpf
We prove the inequality ``$\simgeq$". In the first of \eqref{10.1} it follows from
\begin{equation*}
\begin{split}
\int_{-\infty}^0r^2f^2(\frac t{-r})e^{2r}dr&\geq \int_{-\infty}^{-1}f^2(t)e^{2r}dr
\\&=f^2(t)c,
\end{split}
\end{equation*}
and in the second from
\begin{equation*}
\begin{split}
\int_{-\infty}^0f^2(\frac t{-r})e^{2r}dr&\ge f^2(t)\int_{-1}^0e^{2r}dr+f^2(0)\int_{-\infty}^{-1}e^{2r} dr
\\
&\simgeq f^2(t).
\end{split}
\end{equation*}
We prove the inequality ``$\simleq$". We begin by observing that the conditions $\frac{f}{t}$ and $\frac{f}{t^\eta}$ decreasing are equivalent to  $\frac{f'}{f}\leq \frac1t$ and $\frac{f'}{f}\leq \frac\eta t$ respectively. In the first line is a consequence of 
\begin{equation*}
\begin{split}
\int_{-\infty}^0r^2f^2(\frac t{-r})e^{2r}dr&\simleq \int_{-1}^0f^2(t)e^{2r}dr+\int_{-\infty}^{-1}r^2 f^2(t)e^{2r}dr
\\
&=f^2(t)\Big(\int_{-1}^0 e^{2r}dr+\int_{-\infty}^{-1}r^2e^{2r}dr\Big)
\\
&=f^2(t)c,
\end{split}\end{equation*}
where, the first inequality follows from the fact that $r^2f^2(\frac t{-r})$ on $[-1,0]$ and $f^2(\frac t{-r})$ on $(-\infty,-1]$ achieve their maximum at $-1$ 
(because $s^2f^2(\frac t s)$, $s\in \R^+$, is increasing and $f^2(\frac t s)$ is decreasing) and where, in last equality,  $c$ stands for the sum of the two integrals of the line above.
As for the ``$\simleq$" in the second line, it follows from
\begin{equation*}
\begin{split}
\int_{-\infty}^0f^2(\frac t{-r})e^{2r}dr&=\int_{-1}^0\Big(r^{2\eta} f^2(\frac t{-r})\Big)\frac{e^{2r}}{r^{2\eta}}dr+\int_{-\infty}^{-1}f^2(\frac t{-r})e^{2r}dr
\\
&\leq f^2(t)\Big(\int_{-1}^0 \frac {e^{2r}}{r^{2\eta}}dr+\int_{-\infty}^{-1}e^{2r}dr\Big)
\\
&=:f^2(t)c,
\end{split}\end{equation*}
where the central inequality follows from the fact that  $\Big(r^{2\eta}f^2(\frac t{-r})\Big)$ in  $[-1,0]$ and  $f^2(\frac t{-r})$ in  $(-\infty, -1]$ achieve the maximum at $-1$ (again, because $s^{2\eta}f^2(\frac ts)$ is increasing and $f^2(\frac ts)$ is decreasing); finally, in the last equality, $c$ is a notation for the the sum of the two integrals in the line above.

\epf
\bp
\Label{p10.1}
If $\frac{f}t$ is decreasing, we have
\begin{equation}
\Label{10.2}
\begin{cases}
\no{f(\Lambda)\upiuh}\simleq \no{f(\delta^{-1})\delta\Lambda\upiuh}
\\
\no{f(\delta^{-1})\delta\Lambda\upiuh}\simleq \no{f(\Lambda)u}^b_{-\frac12}.
\end{cases}
\end{equation}
If, instead, $\frac{f}{t^\eta}$ is decreasing and $\eta<\frac12$, then
\begin{equation}
\Label{10.3}
\begin{cases}
\no{f(\Lambda)\upiuh}\simleq \no{f(\delta^{-1})\upiuh},
\\
\no{f(\delta^{-1})\upiuh}\simleq \no{f(\Lambda)u}^b_{-\frac12}.
\end{cases}
\end{equation}
\ep
\bpf
We show how  \eqref{10.2} follows from the first of \eqref{10.1}; (it will be obvious how \eqref{10.3} follows from the second).
We  have
\begin{equation}
\Label{2.7}
\begin{split}
\no{\delta f(\delta^{-1})\Lambda u^{+\,(h)}}^2&=\int_{\R^{2n-1}}\Lambda^2_\xi\int_{-\infty}^0 \delta^2f^2(\delta^{-1})|\tilde u|^2e^{2r\sigma(\dot T)}drd\xi
\\
&=\int_{\R^{2n-1}}\int_{-\infty}^0 s^2\frac{f^2\Big (\frac{\Lambda_\xi}{-s}\Big)}{\Lambda_\xi}e^{2s}|\tilde u|^2d\xi ds
\\
&\sim \int_{\R^{2n-1}}\frac{f^2(\Lambda_\xi)}{\Lambda_\xi}|\tilde u|^2 d\xi,
\end{split}
\end{equation}
where, in the second line, we have used  the substitution $s:=r\sigma(\dot T)$ combined with the fact that $\sigma(\dot T)\sim \Lambda_\xi$ over supp$\,\psi^+$, and, in the third, we have applied Lemma~\ref{l10.1}.
Now, \eqref{2.7}, in the form ``$\simleq$", yields the 
second of \eqref{10.2}. To prove the first of \eqref{10.2}, we observe that
\begin{equation}
\Label{2.8}
\begin{split}
\no{f(\Lambda)u^{+\,(h)}}^2&=\int_{\R^{2n-1}}\int_{-\infty}^0 f^2(\Lambda_\xi)|\tilde u|^2e^{2r\sigma(\dot T)}drd\xi
\\
&\int_{\R^{2n-1}}\int_{-\infty}^0\frac{f^2(\Lambda_\xi)}{\Lambda_\xi}|\tilde u|^2 e^{2s}dsd\xi
\\
&\le \int_{\R^{2n-1}}\frac{f^2(\Lambda_\xi)}{\Lambda_\xi}|\tilde u|^2d\xi,
\end{split}
\end{equation}
where, in the second equality,  we have used the substitution $s:=r\sigma(\dot T)$.
Combination of \eqref{2.7} in the form ``$\simgeq$" with \eqref{2.8} yields the first of \eqref{10.2}.

\epf
\section{Pseudodifferential Calculus with low regularity of coefficients}
\Label{s10}
In this section, $a(x)$ is a function in $H^{1+n-\frac12+\epsilon}$ in $\R^{2n-1}$ which is meant to be the coefficient of a tangential vector field to $b\Om\in H^{2+n +\epsilon}$.  We observe that $a\in H^{1+n-\frac12+\epsilon}$ implies $\widetilde{\Lambda a}\in L^1(\R^{2n-1})$ on account of
\begin{equation}
\begin{split}
\int|\widetilde{\Lambda a}|dx&\le \int\Big|\widetilde{\Lambda a}\Lambda_\xi^{n-\frac12+\epsilon}\Big|\Lambda_\xi^{-(n-\frac12+\epsilon)}d\xi
\\
&\simleq \no{\Lambda a}_{n-\frac12+\epsilon}\no{\Lambda_\xi^{-(n-\frac12+\epsilon)}}_0.
\end{split}
\end{equation}
\bt
\Label{t10.1}
Assume $a\in H^{1+n-\frac12+\epsilon}$ and let $\psi=\psi(\frac\xi{(1+|\xi|^2)^\frac12})$ be a smooth symbol of order $0$.
 Then
\begin{equation}
\begin{split}
\Big|[a(x),\Lambda\Psi]\Big|&\le  c\no{\widetilde{\Lambda a}}_{L^1}
\\
&\simleq c\no{a}_{1+n-\frac12+\epsilon},
\end{split}
\end{equation}
where the inequalities are meant in the sense of operators.
Here, $c:=\sup|\dot\psi|+\sup|\psi|$ and the second inequality follows from the Sobolev immersion Theorem.
\et
\bpf
We start from 
$$
\widetilde{[a(x),\Lambda\Psi]u}=\int\tilde a(\xi-\eta)\Big(\Lambda_\xi\psi(\xi)-\Lambda_\eta\psi(\eta)\Big)\tilde u(\eta)d\eta.
$$
We use the notation $g(\xi):=\Lambda_\xi\psi(\xi)$ and estimate the difference $(g(\xi)-g(\eta))$ which occurs in the formula above. We have 
$$
\di_{\xi_j}g=\frac{\xi_j}{\sqrt{1+|\xi|^2}}\psi+\di_{\xi_j}\psi+\sum_i\di_{\xi_i}\psi\frac{\xi_i\xi_j}{(1+|\xi|^2)}.
$$
Let $v:=\T{grad}\,g$ and note that $|v|$ is bounded in $\R^{2n-1}$; denote by $c$ an upper bound for $|v|$. We have
\begin{equation*}
\begin{split}
g(\xi)-g(\eta)&=\sum_j(\xi_j-\eta_j)\di_{\xi_j}g(\lambda)  
\\
&=:(\xi-\eta)\cdot v(\lambda),
\end{split}
\end{equation*}
where we use the notation   $\lambda=\lambda_{\xi,\eta}$ for the point of mean value. Thus, using Cauchy-Schwarz inequality, we get
\begin{equation*}
\begin{split}
|g(\xi)-g(\eta))|&\le |v| |\xi-\eta|
\\
&\leq c\sqrt{1+|\xi-\eta|^2},
\end{split}
\end{equation*}
for $c= \sup|\dot\psi|+\sup|\psi|$. In conclusion,
\begin{equation*}
\begin{split}
|\widetilde{[a(x),\Lambda\Psi]u}|&\le \int |\tilde a(\xi-\eta)||g(\xi)-g(\eta)||\tilde u(\eta)|d\eta
\\
&\le c\int|\tilde a(\xi-\eta)|\Lambda_{\xi-\eta}|\tilde u(\eta)|d\eta
\\
&\le c\int |\widetilde {\Lambda a}(\xi-\eta)||\tilde u(\eta)|d\eta,
\end{split}
\end{equation*}
and therefore,
\begin{equation*}
\begin{split}
\no{[a(x)\Lambda,\Psi]u}_0&\underset{\T{Plancherel}}=\no{ \widetilde{[a(x)\Lambda,\Psi]u}}_0
\\&\le c\no{|\widetilde{\Lambda a}|*|\tilde u|}\\
&\underset{\T{Young}}\simleq c\no{\widetilde {\Lambda a}}_{L^1}\no{\tilde u}_0
\\
&\underset{\T{Plancherel}}=c\no{\widetilde {\Lambda a}}_{L^1}\no{u}_0.
\end{split}
\end{equation*}

\epf
For any derivative $D=D_{x_i}$, on account of $[D,\psi^+]=0$ we also have 
\begin{equation*}
\begin{split}
[a(x)D,\Psi]&=[a(x),D\Psi]+[D,\Psi a(x)]
\\
&=[a(x),D\Psi]+\Psi Da(x).
\end{split}
\end{equation*}
The main application will be to the commutator 
\begin{equation}
\Label{10.5}
\begin{split}
[\bar L_n,\Psi^+]&=[iT,\Psi^+]
\\
\simleq \no{r}_{2+n +\epsilon},
\end{split}
\end{equation}
where $\Psi^+$ is the operator of positive microlocalization.
\section{The Hardy-Littlewood estimate}
\Label{s20}
Here is the main content of the paper.
\bt
\Label{t2.1}
Let $b\Om\in H^{2+n +\epsilon}$ and suppose $f$, smooth and increasing. Then, for any $u\in C^1(\bar\Om\cap U)$
\begin{equation}\Label{2.3}
\no{f(\Lambda)u}\simleq \no{\delta f(\delta^{-1})\nabla u}+\sum_{j=1}^n\no{\bar L_ju}\quad \T{if $\frac{f}t$ is decreasing,}
\end{equation}
and
\begin{equation}\Label{2.2}
\no{f(\Lambda)u}\simleq \no{f(\delta^{-1})u}+\sum_{j=1}^n\no{\bar L_ju}\quad\T{ if $\frac{f}{t^\eta}$ is decreasing and $\eta<\frac12$}.
\end{equation}
Here, the constants which occur in \eqref{2.3} and \eqref{2.2} are ruled by the $H^{2+n +\epsilon}$ norm of the boundary.
\et
\bpf
We adopt the following terminology. ``Good" is a term controlled by the right side of an estimate and ``neglectable" a term which comes with a small constant or a smaller Sobolev index, or a weaker norm of a term in the left side; this can then be ``absorbed".
We begin by observing that since $|\sigma(T)|\simleq \sum_j|\sigma(\bar L_j)|$ on supp$\,\psi^0$ and $\sigma(T)<0$ on supp$\,\psi^-$, then $\no{u^0}_1+\no{u^-}_1\simleq \sum_j\no{\bar L_ju}+\no{u}$.
(We have to use here the conclusions of Section~\ref{s10} and, namely, \eqref{10.5}.)
 Thus, we only have to prove the theorem for $u^+$. We start by proving \eqref{2.3};   the scheme is to use \eqref{10.2} and to remove the upperfix $(h)$ from both sides. To remove $(h)$ from the left, we have to control by good terms  the difference $u^{+\,(0)}:=u^+-u^{+\,(h)}$.
By the trivial inequality $\no{f(\Lambda)\upiu0}\le\no{\upiu0}_1$, it is enough to control $\no{\upiu0}_1$. We observe that since $\upiu0|_{b\Om}\equiv0$, then the $L_j$'s are controlled by the $\bar L_j$'s, and these are in turn controlled by the single $\bar L_n$ since we are in the microlocal region supp$\psi^+$. This yields the first inequality in the chain below
\begin{equation}
\Label{2.4}
\begin{split}
\no{\upiu0}_1&\simleq\no{\bar L_n\upiu0}+\no{\upiu0}
\\
&\le \no{\bar L_n u^+ }+\no{\bar L_n \upiuh}+\no{u}
\\
&\underset{\T{\eqref{10.5}}}\simleq \T{good}+\no{\bar L_n \upiuh}+\no{u}
\\
&\underset{\T{(i)}}\le \T{good} +\no{(T-\dot T)\upiuh}+c\no{\delta\Lambda\upiuh}
\\
&\underset{\T{(ii)}}\simleq \T{good}+c\no{\delta\Lambda\upiuh}
\\
&\underset{\T{(iii)}}\simleq \T{good}+\no{ u^+ }^b_{-\frac12}
\\
&\underset{\T{(iv)}}\simleq \T{good}+\no{ u^+ }+\no{\Lambda^{-1}\di_r u^+ }
\\
&\underset{\T{(v)}}\simleq \T{good}+\no{ u^+ }+\no{\Lambda^{-1}\bar L_n u^+ }
\\
&=\T{good},
\end{split}
\end{equation}
where (i) follows  from $\bar L_n\upiuh\sim(\bar L_n u)^{+\,(h)}+\tilde \Psi^+$ (for $\Psi^+\prec\tilde \Psi^+$) by Theorem~\ref{t10.1} for a vector field in $\R^{2n-1}$ with parameter $r$ such as $\bar L_n$, (ii) from the fact that $T$ has $C^1$ coefficients, (iii)  from the same change of variables as in Proposition~\ref{p10.1} (cf.  \cite{K02} Lemma~8.4), (iv) is the trace Theorem (cf. \cite{K02} p. 241), and (v) is a consequence of the already mentioned decomposition $\di_r=iT+\bar L_n$. 

We now remove (h) from the right side of the first inequality in \eqref{10.2}. But this is immediate by
\begin{equation*}
\begin{split}
\no{\delta f(\delta^{-1})\Lambda\upiuh}&\le \underset{\T{good}}{\underbrace{\no{\delta f(\delta^{-1})\Lambda u^+}}}+\no{\delta f(\delta^{-1})\Lambda\upiu0}
\\
&\le \T{good}+\no{\upiu0}_1
\\
&\underset{\T{\eqref{2.4}}}\simleq \T{good}.
\end{split}
\end{equation*}

We prove now \eqref{2.2} under the assumption $\frac{f}{t^\eta}$ for $\eta<\frac12$. Recall, again, that this is equivalent to $\frac{f'}{f}\leq\frac{\eta}{t}$. We have
\begin{equation}
\Label{2.12}
\begin{split}
\NO{f(\delta^{-1})\upiuh}&=\left(f(\delta^{-1})\upiuh,f(\delta^{-1})\upiuh\right)
\\
&=\left(\di_r(r)f(\delta^{-1})\upiuh,f(\delta^{-1})\upiuh\right)
\\
&\underset{\T{Integration by parts}}=-2\Re\left(\delta f(\delta^{-1})\upiuh,\di_rf(\delta^{-1})\upiuh\right)
\\
&=-2\Re\left(\delta f(\delta^{-1})\upiuh,(f(\delta^{-1})\di_r+[\di_r, f((\delta)^{-1})])\upiuh\right)
\\
&\simleq2\no{f(\delta^{-1})\upiuh}\Big(\no{\delta f(\delta^{-1})\di_r\upiuh}+\no{\frac{f'}{\delta}\upiuh}\Big).
\end{split}
\end{equation}
Since $\frac{f'}\delta\simleq \eta f$ for $\eta<\frac12$, then the last term between brackets in the bottom of \eqref{2.12} is neglectable and \eqref{2.12} can be rewritten as $\no{f(\delta^{-1})\upiuh}\simleq \no{\delta f(\delta^{-1})\di_r\upiuh}$.
Next,
\begin{equation}
\Label{2.14}
\begin{split}
\no{\delta f(\delta^{-1})\di_r\upiuh}&\leq\no{\delta f(\delta^{-1})T\upiuh}+\no{\delta f(\delta^{-1})\bar L_n\upiuh}
\\
&\simleq\underset{(i)}{\underbrace{\no{\delta f(\delta^{-1})\Lambda u^+}}}+\underset{\T{good by \eqref{2.4}}}{\underbrace{\no{ \upiu0}_1}}+\underset{\T{good}}{\underbrace{\no{\bar L_n u^+}}}.
\end{split}
\end{equation}
To estimate (i), we now apply \cite{K99} Lemma 2.1 to the function $v=f(\delta^{-1})u$ for $\sigma =1,\,s=0$. Note here that 
$v\in C^\infty(\bar \Om\cap U)$ is not needed; what is needed is $|v|\simleq O(\delta^{-\eta})$ for $\eta<\frac12$, which is well fulfilled in our situation. We get
\begin{equation}
\Label{2.15}
\begin{split}
(i)&\simleq \no{f(\delta^{-1}) u^+}+\no{\delta\bar L_n f(\delta^{-1}) u^+}
\\
&=\underset{\T{good}}{\underbrace{\no{f(\delta^{-1}) u^+}}}+\underset{\T{good by \eqref{10.5}}}{\underbrace{\no{\delta f(\delta^{-1})\bar L_n u^+}}}+\no{\delta [\bar L_n,f(\delta^{-1})] u^+}.
\end{split}
\end{equation}
Now,
\begin{equation*}
\begin{split}
|[\delta \di_r,f(\delta^{-1})]|&\le \left|\frac {f'(\delta^{-1})}{\delta}\right|
\\
&\simleq f(\delta^{-1}),
\end{split}
\end{equation*}
and thus the third term in the last line of \eqref{2.15} is also good.
This concludes the proof of Theorem~\ref{t2.1}.

\epf
Assume that  $\frac{f}t$ is decreasing but $\frac{f}{t^\eta}$ is not decreasing for $\eta<\frac12$; then we have \eqref{2.3} but not \eqref{2.2}. However, over harmonic functions these two inequalities are equivalent according to
\begin{equation}
\Label{2.18}
\no{\delta f(\delta^{-1})\nabla u}\sim \no{f(\delta^{-1})u}\quad \T{if $\Delta u=0$.}
\end{equation}
To prove \eqref{2.18}, we use the notation $\di_\nu$ for $\partial r\cdot\nabla$ and remember that $\delta^2 f^2(\delta^{-1})|_{b\Om}\equiv0$, we have
\begin{equation}
\Label{2.5}
\begin{split}
\int_{\Om}\delta^2f^2|\nabla u|^2 dV&=-\int_{\Om}(\delta f \di_\nu u)(f\bar u)dV
+\int_{\Om}\delta^2\frac{f'}{\delta^2}f\di_\nu u\bar u dV+\underset{0}{\underbrace{\int_{\Om}\delta^2f^2\Delta u\bar u dV}}
\\
&\simleq \int_{\Om}(\delta f|\nabla u|)f|\bar u|dV.
\end{split}
\end{equation}
Hence, by Cauchy-Schwarz
$$
\int_\Om\delta^2f^2|\nabla u|^2dV\simleq (\int_\Om \delta^2 f^2 |\nabla u|^2dV)^{\frac12}(\int_\Om f^2|u|^2dV)^{\frac12},
$$
which implies
$$
(\int_\Om \delta^2f^2|\nabla u|^2dV)^{\frac12}\simleq (\int_\Om f^2|u|^2dV)^{\frac12}.
$$

\section{The norm $\no{u}^\Om_f$.}
\Label{snorm}
In these two last sections, $f$ is a smooth increasing function such that $\frac{f}{t^\eta}$ is decreasing for $\eta<\frac{1}{2}$ and hence  $f\simleq t^\eta$, for $t$ large. We pass now from tangential to full pseudodifferential operators. 
 Let $\zeta$ be dual variables to the full system of variables $(x,y)\in\R^{2n}$ and $\Lambda_\zeta=(1+|\zeta|^2)^{\frac12}$, $\zeta^2=\sum_j\zeta_j^2$, be the standard elliptic symbol of order 1 and $\Lambda_{\T{full}}$  the pseudodifferential elliptic operator with symbol $\Lambda_\zeta$. (Here ``full" stresses difference to ``tangential" defined above.) We also consider  some pseudodifferential operator $f(\Lambda_{\T{full}})$ with symbol $f(\Lambda_\zeta)$ for a general  increasing function $f$ defined by $f(\Lambda_{\T{full}})u:=(2\pi)^{-2n}\int_{\R^{2n}}e^{i\langle(x,y),\zeta\rangle}f(\Lambda_\zeta)\tilde u(\zeta)d\zeta$ where $\tilde u$ is the Fourier transform. The space $H^f(\R^n)$ is the space of all $u$ such that $\no{u}_f=\no{f(\Lambda_{\T{full}})u}_0<\infty.$
\bd
\Label{d2.1}
Let $\Om\subset\R^{2n}$;  the space $H^f(\Om)$ is defined as the completion of $C^\infty_c(\Om)$ under the norm $\no{\cdot}_f$. 
\ed


Note that $u\in H^f(\Om)$ does not imply $u|_{b\Om}\equiv0$. To explain this, we consider the characteristic function $\psi$ of the interval $(0,1)$; we have
 $\tilde \psi=\frac{i(e^{-i\xi}-1)}\xi$ and therefore
 \begin{equation}
 \Label{characteristic}
 \begin{split}
 |\tilde\psi|^2(1+|\xi|^{2})^\eta&\simleq(1+|\xi|^2)^{-1+\eta}
 \\
 &\le (1+|\xi|^2)^{-\frac12+\epsilon},
 \end{split}
 \end{equation}
 whose integral is finite. This is similar to the behavior of the Sobolev spaces $H^s(\Om)$ for $s<\frac12$ for which we have $H^s(\Om)=H^s_0(\Om)$ where the suffix $0$ denotes the completion of $C^\infty_c$. However we have introduced a setting which enables to consider  more general  functions such as $f=\log^s(t)$ or $\log^s(\log(t))$.

For a function $u$ in $\Om$, we denote by $\hat u$ the extension by $0$ from $\Om$ to $\R^{2n}$. 
\bp
\Label{pnova}
We have
\begin{equation}
\Label{2.200}
H^f_0(\Om)=\{u:\,\no{\hat u}_f<+\infty\}.
\end{equation}
\ep
\bpf
We prove the inclusion ``$\supset$". Let $\chi_\epsilon\to\delta$ be an approximation of the Dirac measure by $\chi_\epsilon\in C^\infty_c(\B_\epsilon)$ where $\B_\epsilon$ is the $\epsilon$-ball. By a partition of the unity of $\bar \Om$, it is not restrictive to assume that $u\in C^\infty_c(\bar\Om\cap U)$; by taking $U$ sufficiently small, and by  a suitable choice of a ``normal" direction $\nu$, we can assume that $u_\epsilon:=\hat{u}(z+\epsilon\nu)$ satisfies $supp(u_\epsilon)\subset \Omega_{-\frac\epsilon2\nu}:=\{ z\in U:\, r(z+\epsilon\nu)<0\,\}$ so that
$$
\T{supp}(\chi_{\frac\epsilon2}*u_\epsilon)\subset\subset\Om.
$$
It is readily seen that 
$$
\chi_{\frac\epsilon2}*u_\epsilon\to u\quad\T{ in $f$-norm}.
$$
In fact, we pass to Fourier transform and notice that
\begin{equation*}
\widetilde{\chi_{\frac\epsilon2}*u_\epsilon}-\tilde u=(\widetilde{\chi_{\frac\epsilon2}}-1)\tilde u-\widetilde{\chi_{\frac\epsilon2}}(\widetilde{u-u_\epsilon}).
\end{equation*}
Now
\begin{equation*}
\begin{split}
\widetilde{\chi_{\frac\epsilon2}}&=\tilde\chi(\frac\epsilon2\xi)
\\
&=\int e^{i\frac\epsilon2\xi z}\chi(z)dz
\\
&\to \int \chi(z)dz=1\quad\T{pointwise for $\epsilon\to0$},
\end{split}
\end{equation*}
and therefore
$$
\no{(\widetilde{\chi_{\frac\epsilon2}}-1)\tilde u f}_0\to 0\quad\T{by dominated convergence}.
$$
Also,
$$
\sup|\widetilde{\chi_{\frac\epsilon2}}|\le1,\qquad \widetilde{(u_\epsilon-u)}\to0\,\,\T{pointwise},
$$
and therefore
$$
\no{\widetilde{\chi_{\frac\epsilon 2}}(\widetilde{u_\epsilon-u})f}_0\to0\qquad\T{again, by dominated convergence.}
$$

\noindent We prove now ``$\subset$". Now, let $u_m\in C^\infty_c(\Om)$ satisfy $\no{u_m-u}_f\to 0$. This means
\begin{equation*}
\begin{cases}
u_m\underset {L^2}\to \hat u,
\\
f(\Lambda_{\T{full}})u_m\T{ is a Cauchy sequence in $L^2(\R^{2n})$}.
\end{cases}
\end{equation*}
It follows that there exists $v\in L^2(\R^{2n})$ such that $f(\Lambda_{\T{full}})u_m\underset{L^2(\R^{2n})}\to v$; set $w:=\mathcal F^{-1}(\frac {\tilde v}f)$ where $\mathcal F^{-1}$ is the inverse Fourier transform. We have $w=\hat u$; moreover,
 from
 \begin{equation*}
 \begin{split}
 f(\xi)\widetilde{\hat u}&=f(\xi)\tilde w
 \\
 &f(\xi)\frac{\tilde v}{f(\xi)}=\tilde v\in L^2,
 \end{split}
 \end{equation*}
 we conclude that $f(\Lambda_{\T{full}})\hat u\in L^2(\R^{2n})$.
 
 \epf
 Let $\psi_\epsilon$ be the characteristic function of $(0,\epsilon)$; thus $\psi_\epsilon=\psi_1(\frac r\epsilon)$. It is readily seen (cf. \eqref{characteristic}) that $\psi_\epsilon\in H^\eta$ because $\eta<\frac12$ and, by dominated convergence, that $\no{\psi_\epsilon}_f\to0$. 
 Let $\chi_\Om$ and $\chi_{\Om_\epsilon}$ be the characteristic functions of $\Om$ and $\Om_\epsilon:=\{z\in \Om:\,r<-\epsilon\}$ respectively; thus $\hat u=\chi_\Om u$.
 \bp
 \Label{p2.3}
 Let $f$ be increasing and $\frac{f}{t^\eta}$ be decreasing for some $\eta<\frac12$. Let $u\in C^1(\bar\Om)$; then $\no{\chi_{\Om_\epsilon}u-\hat u}\to0$. In particular, $u\in H^f(\Om)$.
 \ep
 \bpf
 We have $\hat u-\chi_\epsilon u=\psi_\epsilon(-r)u$. Clearly,  $\psi_\epsilon\to0$ pointwise; we claim that
 $$
 \no{\psi_\epsilon u}_f<+\infty,
 $$
 which yields $\no{\psi_\epsilon u}_f\to0$ by dominated convergence. To prove the claim we denote by $\tau$ the dual variable  to $r$, write $\eta=\frac12-\alpha$ and note that, for $\tau\ge1$
 \begin{equation*}
 \begin{split}
 |\tau|^{2\eta}\Big|\int_{-\epsilon}^0 e^{-ir\tau}u(r,a)d r \Big|^2&=|\tau|^{2\eta}\Big|\Big[i\frac{e^{-i r \tau}}{\tau} u\Big]^\epsilon_0-\int_0^\epsilon i\frac{e^{-i r \tau}}\tau \di_ r  u(r,a)d r \Big|^2
 \\
 &\simleq |\tau|^{2\eta}\Big(\frac{|e^{-i\epsilon\tau}u(\epsilon,a)-u(0,a)}{\tau^2}+\frac{\sup |\di_ r  u(r,a)|^2}{|\tau|^2}\Big)
 \\
 &\simleq \frac{|\tau|^{2\eta}}{|\tau|^2}\no{u}_{C^1}^2
 \\
 &\simleq |\tau|^{-1-2\alpha}\no{u}^2_{C^1},
 \end{split}
 \end{equation*}
 whose integral is finite.
 \epf
 \bp
 \Label{p2.4}
 Let $f$ be increasing and $\frac{f}{t^\eta}$ be decreasing for some $\eta<\frac12$. Take $u\in D_{\dib}\cap D_{\dib^*}$, the domains of $\dib$ and $\dib^*$ respectively and use the notation  $Q(u,u):=\NO{\dib u}+\NO{\dib^*u}$. Then
 \begin{equation}
 \Label{2.101}
 \NO{u}_f\simleq \NO{f(\Lambda)u}+Q(u,u).
 \end{equation}
 \ep
 \bpf
 
Recall that the regularity of $\Omega$ is at least $C^2$ all through the paper. Also, by Proposition \ref{p2.3} it suffices to prove \eqref{2.101} when $u\in C^1(\bar{\Omega})$; in particular we can use the trace$\backslash$extension argument. Write $\eta=\frac12-\alpha$; we have
 \begin{equation}
 \Label{2.102}
 \frac{f(\Lambda_\zeta)}{\Lambda_\zeta^{\frac12-\alpha}}\le \frac{f(\Lambda_\xi)}{\Lambda_\xi^{\frac12-\alpha}},
 \end{equation}
 because $\Lambda_\zeta\ge \Lambda_\xi$ and $\frac f{t^\eta}$ is decreasing. It follows
 \begin{equation*}
 \begin{split}
 \no{f(\Lambda_{\T{full}})u}&=\no{\frac{f(\Lambda_{\T{full}})}{\Lambda^{\frac12-\alpha}_{\T{full}}}\Lambda^{\frac12-\alpha}_{\T{full}})u}
 \\
 &\underset{\T{\eqref{2.102}}}\simleq \no{\frac{f(\Lambda_{\T{tan}})}{\Lambda^{\frac12-\alpha}_{\T{tan}}}\Lambda^{\frac12-\alpha}_{\T{full}}u}
 \\&\le  \no{\frac{f(\Lambda_{\T{tan}})}{\Lambda_{\T{tan}}^{\frac12-\alpha}}u}_{\frac12-\alpha}^{\T{interpolation}}
 \\&\le \no{\frac{f(\Lambda_{\T{tan}})}{\Lambda_{\T{tan}}^{\frac12}}u}_{\frac12}^{\T{interpolation}}
 \\
 &\underset{\T{\cite{CS00} Lemma~5.1.6}}\le \no{\frac{f(\Lambda_{\T{tan}})}{\Lambda_{\T{tan}}^{\frac12}}u}^b+Q_{\frac {f(\Lambda_{\T{tan}})}{\Lambda^{\frac12}_{\T{tan}}}}(u,u)
 \\
 &\le \no{f(\Lambda_{\T{tan}})u}+\no{\Lambda_{\T{tan}}^{-1}\di_r f(\Lambda_{\T{tan}})u}+Q(u,u)
 \\
 &\simleq\no{f(\Lambda_{\T{tan}})u}+\no{\bar L_n u}+Q(u,u)
 \\
\\&
\simleq\no{f(\Lambda_{\T{tan}})u}+Q(u,u).
\end{split}
\end{equation*}

\epf

\section{An $f(\Lambda)$-estimate for the $\dib$-Neumann problem on a  domain with non-smooth boundary}
\Label{s3}
All through this section, $\Om$ is a bounded pseudoconvex domain of $\C^n$ with $H^{2+n +\epsilon}$  boundary,  $u$ is a form of degree $k\ge1$ and $ f$ a general function with $f$ increasing and $\frac f{t^\eta}$ decreasing  for $\eta<\frac12$ that is, $f'\geq 0$ and $\frac{f'}{f}\leq\frac{\eta}{t}$; in particular $f\lesssim t^\eta$ for $t$ large. Also,  $D_{\dib}$ and $D_{\dib^*}$ denote the domains in $L^2$ of $\dib$ and $\dib^*$ respectively. Finally, let $\no{\cdot}_f$ be the $f$-norm defined in Section~\ref{snorm}.
\bd
{\it $f$-estimate.} We say that an $f$-estimate holds in a neighborhood $U$ of a point $z_o\in b\Om$ when
\begin{equation}
\Label{3.1}
\begin{split}
\NO{u}_f\simleq \NO{\dib u}+\NO{\dib^* u}\quad&\T{for any  form } u\in D_{\dib}\cap D_{\dib^*}
\\
&\T{of degree $k\geq 1$ with supp$\,u\subset\subset U$}.
\end{split}
\end{equation}
\ed
\bd
\Label{d3.2}
{\it $f$-Property. } We say that $\Om$ enjoys the $f$-Property at $z_o$ if there is a family of weights $\phi=\phi_\epsilon\in C^\infty(U\cap \Om)$ for $\epsilon\to0$ which are plurisubharmonic, have uniform bound $|\phi_\epsilon|\le 1$ on $U\cap S_\epsilon$ where $S_\epsilon$ denotes the $\epsilon$-strip $S_\epsilon=\{z\in\Om:\,\delta(z)<\epsilon\}$ and whose Levi-form satisfies
$$
\di\dib \phi_\epsilon\simgeq f^2(\epsilon^{-1})\quad\T{over $ S_\epsilon\cap U$}.
$$
\ed
\bp
\Label{p3.1}
Assume that $\Om$ enjoys the $f$-Property; then there is a single weight $\phi\in C^\infty( \Om\cap U)$ such that 
\begin{equation*}
\begin{cases}
|\phi|\le 1,
\\
\di\dib\phi\ge f_1^2(\delta^{-1}),
\end{cases}
\end{equation*}
where 
\begin{equation}
\Label{3.0}
f_1(t):=\frac {f(t)}{\log(f(t))}.
\end{equation}
\ep
\bpf
Define
$$
\phi=\sum_{k=1}^{+\infty}\frac{\phi_{f^*(2^k)^{-1}}}{\log^22^k},
$$
where $f^*$ is the inverse to $f$.
We first remark that
$\sum^{+\infty}\cdot\simleq \sum^{+\infty}\frac1{k^2}<+\infty$. Moreover, since $\delta \sim f^*(2^k)^{-1}$ and hence $2^k\sim f(\delta^{-1})$ in the difference set $S_{f^*(2^k)^{-1}}\setminus S_{\frac{f^*(2^k)^{-1}}2}$, then
$$
\di\dib\left(\frac{\phi_{f^*(2^k)^{-1}}}{\log^22^k}\right)\simgeq \frac{f^2f^*(2^k)}{\log^22^k}=\frac{f^2(\delta^{-1})}{\log^2(f(\delta^{-1}))}.
$$

\epf

\bt
\Label{t3.1}
Let $\Om\subset\C^n$ be a pseudoconvex, $H^{2+n +\epsilon}$  domain which satisfies the $f$-Property in $U$; then 
we have an $f_1$-estimate for any $u\in D_{\dib}\cap D_{\dib^*}$  with support in $U$.
\et

\bpf
Since $H^{2+n+\epsilon}\subset C^2$, then $C^1(\bar{\Omega})\cap D_{\bar{\partial}}\cap D_{\bar{\partial}^*}$ is dense in the domain of $\bar{\partial}$ and $\bar{\partial}^*$ for the graph norm. Thus, it is sufficient to prove the theorem for the form $u\in C^1_c(\bar{\Omega}\cap U)\cap D_{\bar{\partial}}\cap D_{\bar{\partial}^*}$. Also, under this conditions, we have \eqref{2.101}. If we plug \eqref{2.101} with \eqref{2.3} we get 
\begin{equation}\Label{final} 
 \|u\|_f^2\leq \|f(\delta^{-1})u\|+Q(u,u). 
\end{equation}
By the basic estimate with the weight constructed in Proposition \ref{p3.1}, we can estimate the right hand side of \eqref{final} by $Q(u,u)$. This concludes the proof.  
\epf

\end{document}